\documentclass[12pt]{article}
\usepackage[english]{babel}
\usepackage{amsmath,amsthm}
\usepackage{amsfonts}
\usepackage{graphicx}
\usepackage{enumerate}
\usepackage[usenames,dvipsnames]{color}
\usepackage{subfigure}
\usepackage[right,pagewise,displaymath, mathlines]{lineno}
\usepackage{epstopdf}
\usepackage{color}

\numberwithin{equation}{section}
\numberwithin{figure}{section}
\numberwithin{table}{section}

\newtheorem{theorem}{Theorem}[section]

\newtheorem{definition}[theorem]{Definition}
\newtheorem{example}{Example}[section]

\numberwithin{equation}{section}

\begin{document}

\vspace*{0mm}

\noindent
{\Large \bf Paths and indices of maximal tail dependence}

\vspace*{8mm}

\noindent
{\large  Edward Furman$^{a,*}$, Jianxi Su$^{a}$, and Ri\v{c}ardas Zitikis$^b$}

\bigskip

\noindent
$^a$ Department of Mathematics and Statistics, York University, Toronto, Ontario M3J 1P3, Canada

\noindent
$^b$ Department of Statistical and Actuarial Sciences,
University of Western Ontario, London, Ontario N6A 5B7, Canada

\vspace*{8mm}

\noindent
\rule{165mm}{0.2mm}
\\

\noindent
\textbf{Abstract.}
We demonstrate both analytically and numerically that the existing methods for measuring tail
dependence in copulas may sometimes underestimate the extent of extreme co-movements of
dependent risks and, therefore, may not always comply with the new paradigm of prudent risk
management. This phenomenon holds in the context of both symmetric and
asymmetric copulas with and without singularities. As a remedy, we introduce a notion of \textit{paths of maximal (tail) dependence} and utilize the notion to propose several new indices of tail dependence. The suggested new indices are conservative, conform
with the basic concepts of modern quantitative risk management, and are capable of differentiating between distinct risky positions in situations when the existing indices fail to do so.

\bigskip

%
%

\\

\noindent
\rule{165mm}{0.2mm}

\vfill

\noindent
\rule{25mm}{0.2mm}
\\
{\footnotesize
$^{*}$Corresponding author. Department of Mathematics and Statistics, York University, Toronto, Ontario M3J 1P3, Canada; E-mail: efurman@mathstat.yorku.ca; Tel: 416-736-2100 Ext 33768; Fax: 416-736-5757.
}

\newpage

\section{Introduction}
\label{section-1}

Current regulatory frameworks require enhanced techniques for measuring and managing extremal risks of financial enterprises. In particular, this may involve
\begin{enumerate}[1)]
\item
analyses of marginal risks and dependence structures between them (cf., e.g., McNeil et al., 2005; Jaworski et al., 2010; Jaworski et al., 2013; Durante et al., 2014; and references therein); and
\item
behaviour of aggregated risks, such as sums (cf., e.g., Asimit and Badescu, 2010; Asimit et al., 2011; Embrechts et al., 2013; Embrechts et al., 2014; Puccetti and R\"{u}schendorf, 2014; and references therein).
\end{enumerate}
In many cases, both routes are feasible and complementary to each other. However, there are situations when only one of 1) or 2) is desirable and thus being pursued. We refer to, for example, Biener and Eling (2012), Churchill and Matul (2012), and references therein, for situations in micro-insurance where aggregation is not feasible.

Dependence is quite often modelled using copulas, which have become a well established mathematical tool in actuarial and financial research and practice  (cf., e.g., Frees et al., 1996; Frees and Wang, 2005; McNeil et al., 2005; Embrechts, 2009; Frees, 2010; Jaworski et al., 2010; Jaworski et al., 2013; and references therein), as well as in many other areas such as economics (cf., Patton, 2012), medicine (cf., Nikoloulopoulos and Karlis, 2008), reliability engineering and life sciences (cf., e.g., Balakrishnan and Lai, 2009), and so on. Briefly, a bivariate function $C:[0,1]\times [0,1] \rightarrow [0, 1]$ is a copula if it is grounded, two-increasing, and have uniform marginals (cf., Nelsen, 2006).

When quantifying co-movements of extreme risks, the behaviour of $C$ around the upper-right and lower-left vertices of the domain of definition $[0,1]\times [0,1]$ gives rise to indices of upper and lower tail dependence. While these indices are concerned with phenomena in opposing risk-tails, they can be formally unified by shifting attention from the underlying copula to the corresponding survival copula.
For example, given the index of lower-tail dependence $\lambda_L:=\lambda_L(C)$, which is
\begin{equation}\label{ind-ll}
\lambda_L=  \lim_{u \downarrow 0} {C(u,u)\over u},
\end{equation}
the corresponding index of upper-tail dependence $\lambda_U:=\lambda_U(C)$ is
\[
\lambda_U= \lim_{u \downarrow 0} {1-2(1-u)+C(1-u,1-u)\over u}
= \lim_{u \downarrow 0} {\widehat{C}(u,u)\over u} ,
\]
where $\widehat{C}$ is the survival copula of $C$; hence, $\lambda_U(C)=\lambda_L(\widehat{C})$. (Throughout the paper we use `$:=$' instead of the customary equality sign `$=$' when we want to emphasize that equality is \textit{`by definition.'})

Similar arguments hold for other indices of tail dependence, such as the index of weak lower-tail dependence $\chi_{L}:=\chi_{L}(C)$ given by
\begin{equation}\label{ind-cl}
\chi_{L}=\lim_{u\downarrow 0}\frac{2\log u }{\log C(u,u)}-1 ,
\end{equation}
and its upper variant $\chi_{U}:=\chi_{U}(C)=\chi_{L}(\widehat{C})$ (cf., Coles et al., 1999; Heffernan, 2000; Fischer and Klein, 2007; and references therein).
In view of these notes, we restrict our following considerations to the behaviour of copulas near the lower-left vertex of their domain of definition.

Instead of comparing $C(u,u)$ and $u \in (0,1)$ as in limits (\ref{ind-ll}) and  (\ref{ind-cl}), we can more generally explore the asymptotic behaviour of $C(u,u)$ in the form
\begin{equation}\label{ind-kl}
C(u,u)= \ell(u) u^{\kappa_L} \quad \textrm{when} \quad u\downarrow 0,
\end{equation}
which defines the index $\kappa_L:=\kappa_L(C)\in [1,\infty )$ of lower-tail dependence, assuming of course that equation (\ref{ind-kl}) holds for a slowly varying at $0$ function $\ell(u)$ (cf., Ledford and Tawn, 1996). We note that smaller values of
$\kappa_L$ correspond to stronger interdependences, thus motivating the search for the smallest possible analogues of $\kappa_L$, and this is our main goal in the present paper.

The rest of the paper is organized as follows. In Section \ref{section-2} we describe the main framework for our ideas, including the class of admissible dependence paths and its subclass of maximal dependence paths. In Section \ref{num-ill} we  present a numerical example that shows that the existing indices of tail dependence may underestimate the extent of extreme co-movements of dependent risks. In Section \ref{sub-mo} we explore several families of copulas for which paths of maximal dependence can be derived in closed form, whereas somewhat more intricate examples are given in Section \ref{fe-2}. In Section \ref{section-4}, we further extend our general framework by introducing new tail orderings of distinct copulas.  Section \ref{section-7}  concludes the paper. Tedious calculations are relegated to Appendix \ref{appendix}.

\section{Paths of maximal dependence and related indices}
\label{section-2}

In what follows we work with bivariate copulas $C:[0,1]\times [0,1] \rightarrow [0, 1]$ only, because they convey the main idea of the present paper in a simple and illuminating way. With the notation $R(u,v)=[0,u]\times [0,v]$, the copula-value $C(u,v)$ is the probability that the bivariate random vector $(U,V)$ with uniform marginals $U$ and $V$ falls into the rectangle $R(u,v)$, that is,
\begin{equation}\label{main-0}
C(u,v)=\mathbf{P}\big [(U,V) \in R(u,v) \big ].
\end{equation}
The class of all rectangles $R(u,v)$ contains the subclass of all squares $R(u,u)$, and these are the ones that have traditionally been used for measuring the strength of tail dependence. Namely, the classical indices of lower-tail dependence are based on the behaviour of probability (\ref{main-0}) when the rectangle $R(u,v)$ shrinks along the diagonal
\begin{equation}\label{diagonal}
\{(u,v) \in [0,1]^2 : u=v\}
\end{equation}
when $u \downarrow 0$. (Throughout the paper, we work only with diagonal (\ref{diagonal}) and thus call it \textit{the} diagonal.) In the context of the present paper, the diagonal is a path of (tail) dependence, but there
are of course many other possible paths and we shall next describe them.

In order to make the problem well posed, some restrictions on the class of possible paths must be imposed. First, we observe that for the independence copula, it is natural to require that every path would reflect the same degree of tail dependence, that is, would have the same probability (\ref{main-0}). For two functions $\varphi, \psi:[0,1]\rightarrow [0,1]$, this implies that every path $(\varphi(u), \psi(u))_{0\le u \le 1}$ has to be necessarily of the form  $(\varphi(u),u^2/\varphi(u))_{0\le u \le 1}$, that is, we must have  $\psi(u)=u^2/\varphi(u)$. Certainly, both $\varphi(u)$ and $u^2/\varphi(u)$ must be in the interval $[0,1]$, and thus $\varphi(u)$ must  be in the interval $[u^2,1]$, which justifies the first property of the following definition.

\begin{definition}\label{def-ap}\rm
We call a function $\varphi: [0,1] \to [0,1]$ admissible if it satisfies the following two properties:
\begin{enumerate}[(1)]
\item
$\varphi(u) \in [u^2,1] $ for every $u\in [0,1] $; and
\item \label{condition-2}
both $\varphi(u)$ and $u^2/\varphi(u) $ converge to $0$ when $u\downarrow 0$.
\end{enumerate}
We call the path $(\varphi(u),u^2/\varphi(u))_{0\le u \le 1}$ admissible whenever the function $\varphi $ is admissible. Finally, we use $\mathcal{A} $ to denote the set of all admissible functions $\varphi $.
\end{definition}

The second property of Definition \ref{def-ap} is related to the fact that we are interested in the behaviour of copulas near the lower-left vertex of their domain of definition.

Motivated by the idea of determining the strongest extreme co-movements of risks, among all admissible functions $\varphi \in \mathcal{A}$, we search for those that maximize the probability
\begin{equation}
\Pi_{\varphi }(u) =C\big (\varphi(u),u^2/\varphi(u)\big )
\label{probab}
\end{equation}
or, equivalently, the distance function
\begin{equation}
\label{distf}
d_\varphi\left(
C, C^\perp
\right)(u)=
C\big (\varphi(u),u^2/\varphi(u)\big )-C^{\perp}(\varphi(u),u^2/\varphi(u)),
\end{equation}
where $C^\perp$ is the independence copula, i.e., $C^\perp(u, v)=uv$
for all $0\leq u,v\leq 1$. Obviously, the function $\varphi_0(u)=u $ is  admissible and yields
the representation of the diagonal path that serves as a building block for the classical indices of lower-tail dependence. However, this path may not maximize probability (\ref{probab}), as we illustrate in following Example \ref{ExMO1}. In fact, in view of (\ref{distf}), the classical index $\kappa_L$ may serve as neither a maximal nor minimal measure of tail dependence of the copula $C$.

\begin{example}
\label{ExMO1}
\rm
Consider the Marshall-Olkin copula
\begin{equation}\label{copula-mo}
C_{a,b}(u,v)=\min(u^{1-a}v,uv^{1-b})
\quad \textrm{for} \quad 0\le u,v\le 1,
\end{equation}
where $a,b\in [0,1]$ are parameters. We check that $C_{a,b}(u,u)=\min(u^{2-a},u^{2-b})$ and thus
\[
\Pi_{\varphi_0 }(u)=u^{\kappa_L} \quad \textrm{with} \quad \kappa_L=2-\min\{a,b\}.
\]
Next we choose the admissible function $\varphi_1(u)=u^{2b/(a+b)}$. See Figure \ref{fig-mo} for the corresponding path,
\begin{figure}[t!]
\centering
\includegraphics[width=7.5cm,height=7.5cm]{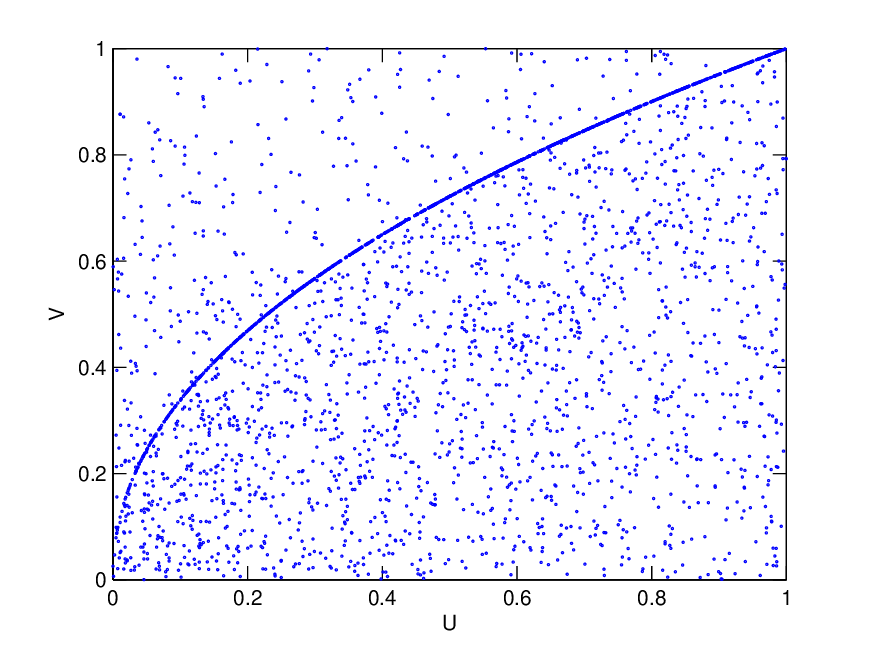}
\hfill
\includegraphics[width=7.5cm,height=7.5cm]{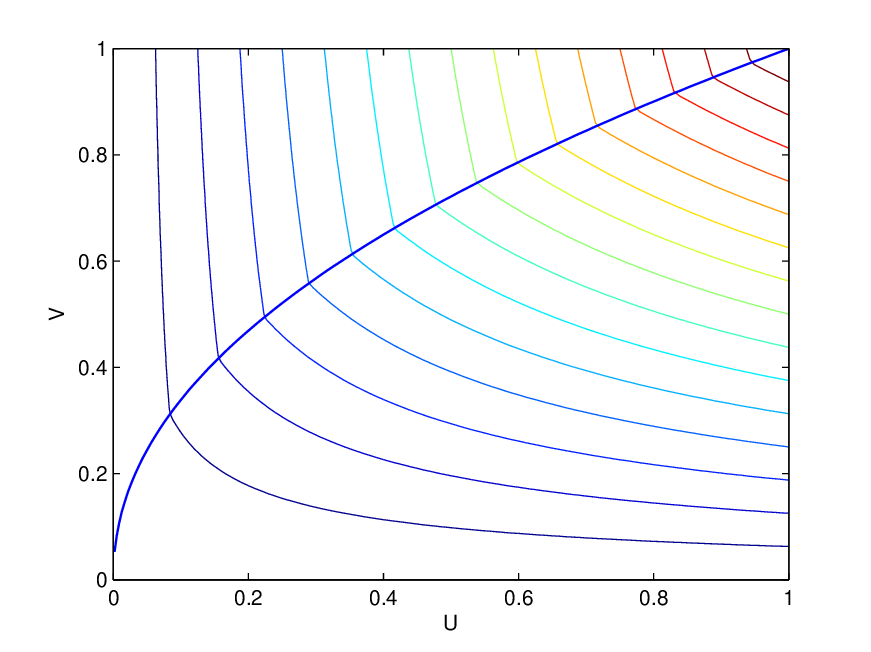}
\caption{Simulated Marshall-Olkin copula (left) and its contours (right) when $a=0.3529$ and $b=0.75$, with the path of maximal dependence  $(\varphi_1(u),u^2/\varphi_1(u))_{0\le u \le 1}$ superimposed on the right-hand panel.}
\label{fig-mo}
\end{figure}
which we later show to be the (only) path of maximal dependence. We check that
\begin{equation}\label{mo-1}
\Pi_{\varphi_1 }(u)=u^{\kappa_L^*} \quad \textrm{with} \quad
\kappa_L^*=2-{2ab \over a+b } .
\end{equation}
Clearly, $\kappa_L^*\le \kappa_L$ with the equality holding only when $a=b$, that is, when the Marshall-Olkin copula is symmetric.
\end{example}

This example  motivates the following definition of a subclass of admissible functions.

\begin{definition}\label{def-ltd}\rm
Given a copula $C$, an admissible function $\varphi^* \in \mathcal{A} $
is called a function of maximal dependence if
\[
\Pi_{\varphi^* }(u)=\max_{\varphi \in \mathcal{A} }\Pi_{\varphi }(u)
\quad \textrm{for all} \quad u\in [0,1].
\]
We conveniently use the simpler notation $\Pi^*(u)$ instead of $\Pi_{\varphi^* }(u)$, but may enhance it to $\Pi^*(u \mid C)$ when a need arises to emphasize the dependence of $\Pi^*(u)$ on $C$. Also, we refer to $(\varphi^\ast(u), u^2/\varphi^\ast(u))_{0\leq u\leq 1}$ as a path of maximal dependence.
\end{definition}

A useful technique for deriving function(s) and thus path(s) of maximal dependence -- and we have frequently used it in our explorations -- is based on:
\begin{itemize}
\item
searching for critical points of the function $x\mapsto C(x,u^2/x)$ over the interval $[u^2,1]$, and for each $u\in [0,1]$, and then
\item checking which of the solutions is/are global maximum/maxima.
\end{itemize}
Accomplishing these tasks may sometimes be relatively easy (Section \ref{sub-mo}), sometimes challenging (Furman et al., 2015), and in some cases obtaining closed-form solutions may not even be possible (Section \ref{fe-2}). Sometimes, especially when formulas for conditional copulas are readily available, it is useful to recall that partial derivatives of copulas are conditional copulas, and thus the task of determining the set of critical points becomes equivalent to finding all the solutions in $x\in [u^2,1]$ to the equation
\begin{equation}\label{con-1}
    xC_{2\mid 1}\left ({u^2\over x}  \mid x \right )
    = {u^2\over x}C_{1\mid 2}\left (x \mid {u^2\over x} \right ).
\end{equation}
For example, equation (\ref{con-1}) has played a pivotal role when handling the Gaussian copula by Furman et al. (2015).

{\color{blue}
Functions of maximal dependence may or may not be continuous, and this is related to uniqueness of the path of maximal dependence. To demonstrate this fact, we start with a theorem whose proof is given in Appendix \ref{appendix}.

\begin{theorem}\label{cont}
If a copula has only one function of maximal dependence, then the function is continuous.
\end{theorem}

}

For copulas with at least two paths of maximal dependence, there are automatically infinitely many paths of maximal dependence, as has been very rightly noted by one of the referees of the current paper: given two functions $\varphi^*_1 $ and $\varphi^*_2 $ of maximal dependence, their combination $\varphi^*_A(u):=\varphi^*_1 (u)\mathbf{1}_{A}(u)+\varphi^*_2(u)\mathbf{1}_{A^c}(u) $ is a function of maximal dependence for any $A\subset [0,1]$, where $\mathbf{1}_{A}$ and $\mathbf{1}_{A^c}$ are the indicator functions of the set $A$ and its complement $A^c$, respectively. Of course, functions such as $\varphi^*_A$ are perfectly legitimate and the probability of maximal dependence will be the same irrespective of the chosen path: $\Pi_{\varphi^*_A }=\Pi_{\varphi^*_1 }=\Pi_{\varphi^*_2 }$. Consequently, the new indices of tail dependence that we introduce next do not depend on the chosen path of maximal dependence.

We now use Definition \ref{def-ltd} to introduce
three conservative variants of  the classical indices $\lambda_L$, $\chi_L$ and $\kappa_L$. Namely, assuming that the limits below exist, we set
\begin{equation}
\label{lambdamax}
\lambda_L^*:=\lambda_L^*(C)= \lim_{u \downarrow 0} {\Pi^*(u)\over u},
\end{equation}
\[
\chi_L^*:=\chi_L^*(C)=\lim_{u\downarrow 0}\frac{2\log u }{\log \Pi^*(u)}-1,
\]
and let $\kappa_L^*:=\kappa_L^*(C)$ be such that
\begin{equation}\label{mo-1b}
\Pi^*(u)=\ell^*(u) u^{\kappa_L^*}\quad \textrm{when} \quad u\downarrow 0
\end{equation}
for a slowly varying at $0$ function $\ell^* (u)$, assuming that such a function exists. In our illustrative examples below, we  concentrate on calculating $\kappa_L^*$ only.

We conclude this section by noting that this paper is not the first attempt to abandon the diagonal section of copulas when measuring tail dependence.
For example, in their research of large claims reinsurance, Asimit and Jones (2008) rely on the asymptotic behaviour of copulas along non-diagonal paths. Asimit and Badescu (2010) rely on the extreme behaviour of copulas along a variety of paths when exploring a time dependent risk model with dependent inter-claim times and claim amounts. Joe et al. (2010) introduce the tail dependence function $b(w_1, w_2;\ C)= \lim_{u\downarrow 0}C(uw_1,uw_2)/u $ for $ w_1>0$ and $w_2>0$. We refer to Asimit et al. (2011),  Weng and Zhang (2012), and Li and Wu (2013) for further developments on the topic. Hua and Joe (2014) use the excess-of-loss economic pricing functional to propose and study a measure of tail dependence that does not rely on the diagonal dependence path. There are also several other related works but none of them -- due to different research goals -- aim at maximal-dependence paths and, in turn, at
indices of maximal dependence.

\section{Numerical illustration}
\label{num-ill}

Here we present a numerical example that questions the decisive role of the diagonal path in measuring tail dependence in copulas.  We nevertheless stress at the outset that our discussion is not a criticism of the role of the diagonal when investigating copulas -- it does play a pivotal role in the analysis of a variety of other aspects as elucidated by, e.g., Durante et al. (2014) and references therein.

Let $X$ and $Y$ be two random variables, which we assume to follow the same Pareto-II (also known as Lomax)  distribution. Its  decumulative distribution and probability density functions are given by, respectively,
\begin{equation}\label{1d-pareto}
\bar{F}(x)  =\left(  \dfrac{x-\mu}{\sigma}+1\right)  ^{-\alpha}
\quad \textrm{and} \quad
f(x)  =\dfrac{\alpha}{\sigma}\left(  \dfrac{x-\mu}{\sigma}+1\right)  ^{-\left(  \alpha+1\right)  }
\end{equation}
for all $x\geq\mu\in(-\infty, \infty)$. Throughout this section, we set $\mu=0$ and $\sigma =1$ for simplicity, and set $\alpha =4$ to ensure the finiteness of all quantities that we consider.

Let the dependence structure between $X$ and $Y$ be given by the Marshall-Olkin copula; see equation (\ref{copula-mo}). Hence, the joint cumulative distribution function of $X$ and $Y$ is equal to $C_{a,b}(F(x),F(y)), 0\leq x,y< \infty$ (for applications of this model to insurance, we refer to, e.g., Asimit et al., 2010). Following Embrechts et al. (2003), we set $a=0.3529$.  For several values of $b\in [0,1]$, we calculate
\begin{enumerate}[(i)]
\item \label{k}
Kendall's $\tau =\tau(C_{a,b})$ index of dependence
\item
classical index $\kappa_L =\kappa_L(C_{a,b})$ of lower-tail dependence given by equation (\ref{ind-kl})
\item
newly suggested index $\kappa_L^* =\kappa_L^*(C_{a,b})$ defined by equations (\ref{mo-1}) and (\ref{mo-1b})
\item
value-at-risk
$\mathrm{VaR}_q[Z]= \inf\{x\in\mathbf{R}:F_Z(x)\geq q\}$ for $Z=X+Y$
\item
conditional tail expectation
$\mathrm{CTE}_q[Z]= \mathbf{E}[Z|\ Z>\mathrm{VaR}_q[Z]]$
\item \label{mtv}
modified tail variance
\[
\mathrm{MTVar}_q[Z]= \mathrm{CTE}_q[Z]+\frac{1}{\mathrm{CTE}_q[Z]}\mathbf{Var}[Z|\ Z>\mathrm{VaR}_q[Z]]
\]
\end{enumerate}
where $q\in(0,1)$. When conducting our numerical calculations, we set $q$ to $0.990$ and $0.995$ in all of the aforementioned weighted
risk measures (cf., e.g., Furman and Zitikis, 2010; and references therein). We have summarized the results in Table \ref{tbl-1}.
\begin{table}[h!]
\centering
\caption{Quantities (\ref{k})--(\ref{mtv}) for the Marshall-Olkin copula with Pareto-II marginals.}\bigskip
\begin{tabular}{c|c|c|c|c|c|c|c}
\hline\hline
\multicolumn{2}{c|}{Parameters} & \multicolumn{3}{c|}{Indices of dependence} & \multicolumn{3}{|c}{
Risk measures} \\
$q$ & $b$ & $\tau$ & $\kappa_L$ & $\kappa_L^*$& $\mathrm{VaR}_q[Z]$ & $\mathrm{CTE}_q[Z]$ & $\mathrm{MTVar}_q[Z]$ \\
\hline\hline\
  $0.9900$ & $0.7500$ & $0.3158$ & $1.6471$ & $1.5200$ & $3.4621$ & $4.8599$ & $5.5808$ \\
   & $0.5000$  & $0.2609$ & $1.6471$ & $1.5862$ &	$3.4095$ & $4.7606$ & $5.4691$ \\
   & $0.3529$ & $0.2143$ & $1.6471$ & $1.6471$ & $3.3612$ & $4.6926$ & $5.3951$ \\
 \hline
 $0.9950$ & $0.7500$ & $0.3158$ & $1.6471$ & $1.5200$ & $4.2925$ & $5.8976$ & $6.7004$ \\
 & $0.5000$ & $0.2609$ & $1.6471$ & $1.5862$ & $4.2114$ & $5.7782$ & $6.5552$ \\
 & $0.3529$ & $0.2143$ & $1.6471$ & $1.6471$ & $4.1460$ & $5.6801$ & $6.4268$ \\
 \hline\hline
\end{tabular}
\label{tbl-1}
\end{table}
\noindent

Noting that smaller values of $\kappa_L^*$ mean stronger tail dependence, it is illuminating to observe from Table \ref{tbl-1} that the smaller the values of $\kappa_L^*$ are, the larger the values of VaR, CTE and MTVar are, whereas the classical index $\kappa_L$ does not change.

\section{Examples}
\label{sub-mo}

We start with several families of copulas for which paths of maximal dependence are derivable in closed form and with moderate amount of effort.

\subsection{Marshall-Olkin copula}

Recall that the Marshall-Olkin copula is defined by formula (\ref{copula-mo}). Next, for every $u\in [0,1]$, the function $x\mapsto C_{a,b}(x,u^2/x)$ defined on the interval $[u^2,1]$ is equal to $u^{2(1-b)}x^b$ for all $x\le x_0=u^{2b/(a+b)}$ and $u^2/x^a$ for all $x\ge x_0$. Hence, the unique maximum of the function is achieved at the point $x=x_0$, and thus the function of maximal dependence of the Marshall-Olkin copula is unique and given by
\begin{equation}
\varphi^*(u)=u^{2b/(a+b)}.
\label{qq-3omo0}
\end{equation}
Consequently, the maximal probability is
\begin{equation}
\Pi^*(u)=u^{2-2ab/(a+b) },
\label{qq-3omo1}
\end{equation}
and thus the lower-tail index of maximal dependence  is
\begin{equation}
\kappa^*_L =2- {2ab \over a+b}.
\label{qq-3omo}
\end{equation}

\subsection{Mixture of Marshall-Olkin copulas}

\begin{figure}[t!]
\centering
\includegraphics[width=7.5cm,height=7.5cm]{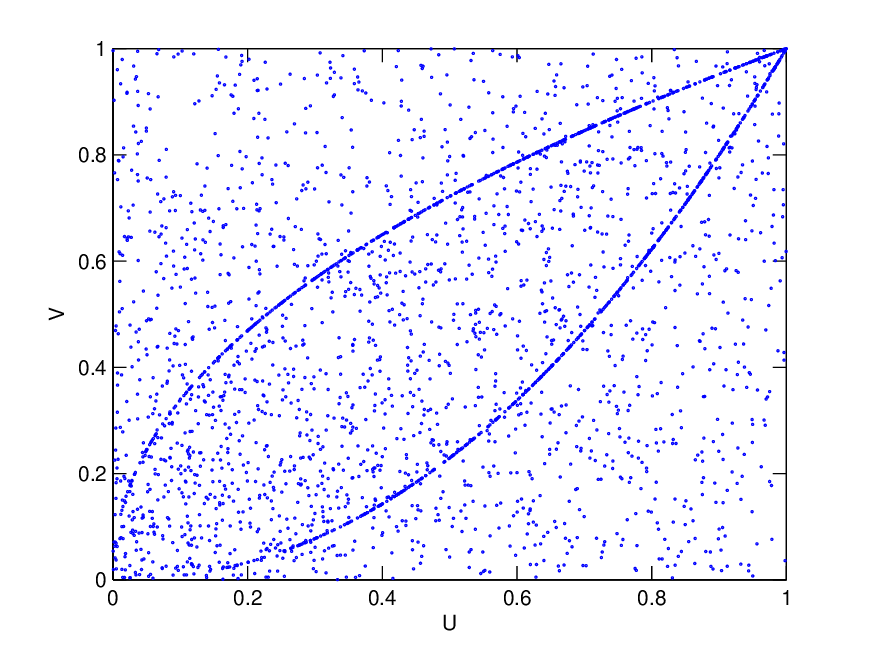}
\hfill
\includegraphics[width=7.5cm,height=7.5cm]{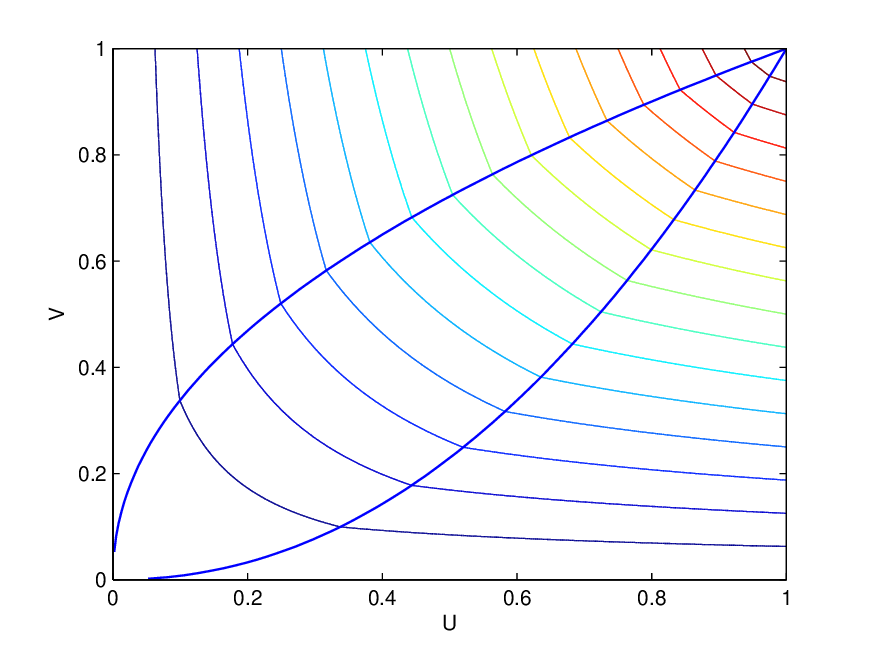}
\caption{Simulated mixture of two Marshall-Olkin copulas (left) and its contours (right) when $a=0.3529$ and $b=0.75$, with the two paths of maximal dependence $(\varphi_k^*(u),u^2/\varphi_k^*(u))_{0\le u \le 1}$,  $k=1,2$, superimposed on the right-hand panel.}
\label{fig-mixmo}
\end{figure}

We see from formula (\ref{qq-3omo0}) that the path of maximal dependence is diagonal if and only if $a=b$, and thus if and only if the Marshall-Olkin copula is symmetric. We  next show that this fact cannot be generalized to arbitrary symmetric copulas. Namely, there are symmetric copulas whose paths of maximal dependence are not diagonal. To show this, we use the $0.5/0.5$ mixture of two mirrored (around the diagonal) Marshall-Olkin copulas, that is,
\begin{equation}\label{mo-1aa0}
C(u,v)= {1\over 2} \Big ( C_{a,b}(u,v)+C_{b,a}(u,v) \Big )
\quad \textrm{for} \quad 0\le u,v\le 1,
\end{equation}
with some $a\ne b$. The copula $C(u,v)$ is symmetric. Tedious calculations (Appendix \ref{appendix}) show that there are two paths of maximal dependence (Figure \ref{fig-mixmo})
given by the functions
\begin{equation}\label{mo-1a0}
\varphi_{1}^*(u)=u^{2b/(a+b)}  \quad \textrm{and} \quad \varphi_{2}^*(u)=u^{2a/(a+b)}.
\end{equation}
Hence, none of the paths of maximal dependence coincides with the diagonal because $a\ne b$, and for any of the two paths, the maximal probability is
\begin{align}
\Pi^*(u)
&={1\over 2}  u^{2-2ab/(a+b)} \Big ( 1+ u^{2(ab-\min\{a,b\})/(a+b)} \Big )
\notag
\\
&= {1\over 2}  u^{2-2ab/(a+b)} \big ( 1+ o(1) \big ) \quad \textrm{when} \quad u\downarrow 0.
\label{qq-3mixmo1}
\end{align}
Hence, the lower-tail index of maximal dependence is
\begin{equation}
\kappa^*_L =2- {2ab \over a+b},
\label{qq-3eee}
\end{equation}
which is the same as in equation (\ref{qq-3omo}) for the ordinary Marshall-Olkin copula.

\subsection{Farlie-Gumbel-Morgenstern copula}

We saw in the previous subsection that the symmetry of a copula does not necessarily imply that the diagonal is a path of maximal dependence. Of course, some symmetric copulas do have diagonal paths of maximal dependence. One of such examples is the already discussed Marshall-Olkin copula when $a=b$. Another  example is the Farlie-Gumbel-Morgenstern (FGM) copula
\[
C(u,v)=uv\big ( 1+\alpha(1-u)(1-v) \big )
\quad \textrm{for} \quad 0\le u,v\le 1,
\]
with parameter $\alpha \in [-1,1]$. Our choice to include the FGM copula among illustrative examples has been motivated by the fact that within the context of the present paper, it sometimes matters whether a copula is or is not positively quadrant dependent (PQD). Specifically for the FGM copula, when $\alpha<0$, then the copula is negatively quadrant dependent (NQD); when $\alpha>0$, then it is  PQD; and when $\alpha=0$, then it is the independence copula. We shall see the role of the sign of $\alpha $ in our following considerations.

To find the path of maximal dependence for the FGM copula, for every $u\in [0,1]$, we first search for those $x\in [u^2,1]$ that solve the equation $(\partial /\partial x )C(x, u^2/x)=0$. Since
\begin{align*}
\frac{\partial C(x, u^2/x)}{\partial x}
&=\frac{\partial}{\partial x} u^2\big ( 1+\alpha (1-u^2/x)(1-x)\big )
\\
&=u^2 \alpha(-1+u^2/x^2),
\end{align*}
we have $(\partial /\partial x )C(x, u^2/x)=0$ if and only if $x=u$. Therefore, when $\alpha >0$, then the diagonal is the unique path of maximal dependence, that is, $\varphi^*(u)=u$. In this case, the maximal probability is $\Pi^*(u)=u^2( 1+\alpha(1-u)^2) $, and thus the lower-tail index of maximal dependence is
\[
\kappa^*_L=2.
\]
As to the case $\alpha \le 0$, we first note that when $\alpha=0$, then every admissible path is a path of maximal dependence. When $\alpha<0$, then $C(x,u^2/x)$ reaches its maximal value at either $u^2$ or $1$ or both, but neither of the functions $\varphi(u)=u^2$ and $\varphi(u)=1$ is admissible because they fail to satisfy property (2) of Definition \ref{def-ap}.

\section{Lower-tail dependence comparisons}
\label{section-4}

We  now look at the index $\lambda_L^*$ (cf. equation (\ref{lambdamax})) from a different angle. To this end we first note that for the Fr\'{e}chet upper bound copula $C^\top(u,v)=\min\{u,v\}$,
$0\leq u,v\leq 1$, the path of maximal dependence is diagonal. Consequently, $\Pi^*(u \mid C^\top )=u$ and so, given any copula $C$, the index $\lambda_L^*:=\lambda_L^*(C)$ can be rewritten as
\[
\lambda_L^*(C)=\lim_{u \downarrow 0}{\Pi^*(u \mid C)\over \Pi^*(u \mid C^\top)}.
\]
This suggests that, given two copulas $C_1$ and $C_2$, we can compare their lower-tail dependencies using the index
\[
\lambda_L^*(C_1,C_2)= \lim_{u \downarrow 0} {\Pi^*(u \mid C_1) \over \Pi^*(u \mid C_2) }  .
\]
Obviously, $\lambda_L^*(C)=\lambda_L^*(C,C^\top)$.

\begin{definition}
The copula $C_1$ is said to be more lower-tail maximally-dependent (LTMD) than the copula $C_2$ if $\lambda_L^*(C_1,C_2)>1$; and less LTMD if $\lambda_L^*(C_1,C_2)<1$.
\end{definition}

An analogous interpretation can be made about the index $\chi_L^*:=\chi_L^*(C)$, but now using the independence copula $C^\perp(u,v)=uv, 0\leq u,v\leq 1$, for which every admissible path is a path of maximal dependence. Hence,
\[
\chi_L^*(C)=\lim_{u\downarrow 0}\frac{\log\Pi^*(u \mid C^\perp)}{\log \Pi^*(u \mid C)}-1.
\]
We generalize this index to the case of two arbitrary copulas $C_1$ and $C_2$ by introducing the index
\[
\chi_L^*(C_1,C_2)= \lim_{u \downarrow 0} {\log \Pi^*(u \mid C_2) \over \log \Pi^*(u \mid C_1) }  -1.
\]
Obviously, $\chi_L^*(C)=\chi_L^*(C,C^\perp)$.

\begin{definition}
The copula $C_1$ is said to be more weakly lower-tail maximally-dependent (WLTMD) than the copula $C_2$ if $\chi_L^*(C_1,C_2)>0$; and less WLTMD if $\chi_L^*(C_1,C_2)<0$.
\end{definition}

We can now compare all the specific copulas considered in this paper, as well as other ones. Without getting too specific at the moment, assume that $C_1$ and $C_2$ are such that, for some constants $c^*(C_k)$ and $\kappa_L^*(C_k)$, $k=1,2$, the asymptotic formulas $\Pi^*(u\mid C_1) = c^*(C_1) u^{\kappa_L^*(C_1)}(1+o(u))$ and $\Pi^*(u\mid C_2) = c^*(C_2)u^{\kappa_L^*(C_2)}(1+o(u))$ hold when  $u \downarrow 0 $. Two cases follow:

\begin{itemize}
\item
When $\kappa_L^*(C_1)=\kappa_L^*(C_2)$, then
\[
\lambda_L^*(C_1,C_2)= {c^*(C_1) \over c^*(C_2) }.
\]
To illustrate, we note that the classical Marshall-Olkin copula, which we denote here by $C_{\textsc{mo}}$, and mixture (\ref{mo-1aa0}) of two Marshall-Olkin copulas, which we denote here by $C_{\textsc{mixmo}}$, have the same $\kappa_L^*$-indices. Equations (\ref{qq-3omo1}) and (\ref{qq-3mixmo1}) imply that $\lambda_L^*(C_{\textsc{mo}},C_{\textsc{mixmo}})=2$, which we interpret as saying that $C_{\textsc{mo}}$ is more LTMD than $C_{\textsc{mixmo}}$ whenever $a\neq b$. When $a=b$, then we of course have $\lambda_L^*(C_{\textsc{mo}},C_{\textsc{mixmo}})=1$ because $C_{\textsc{mixmo}}=C_{\textsc{mo}}$ in this case.
\item
When $\kappa_L^*(C_1)\ne \kappa_L^*(C_2)$, then only the index $\chi_L^*(C_1,C_2)$ is of interest, which can be expressed by the formula
\[
\chi_L^*(C_1,C_2)= {\kappa_L^*(C_2) \over \kappa_L^*(C_1) }  -1.
\]
In particular, we have $\chi_L^*(C)=\chi_L^*(C,C^\perp)=2/\kappa_L^*(C) -1$.
\end{itemize}

\section{Further examples}
\label{fe-2}

For all the copulas hitherto, we have derived closed-form expressions for their paths of maximal dependence. This may not always be the case as we illustrate with following examples; yet, we shall be able to derive closed-form expressions for the index $\kappa_L^*$.

\subsection{Generalized Clayton copula}

\begin{figure}[t!]
\centering
\includegraphics[width=7.5cm,height=7.5cm]{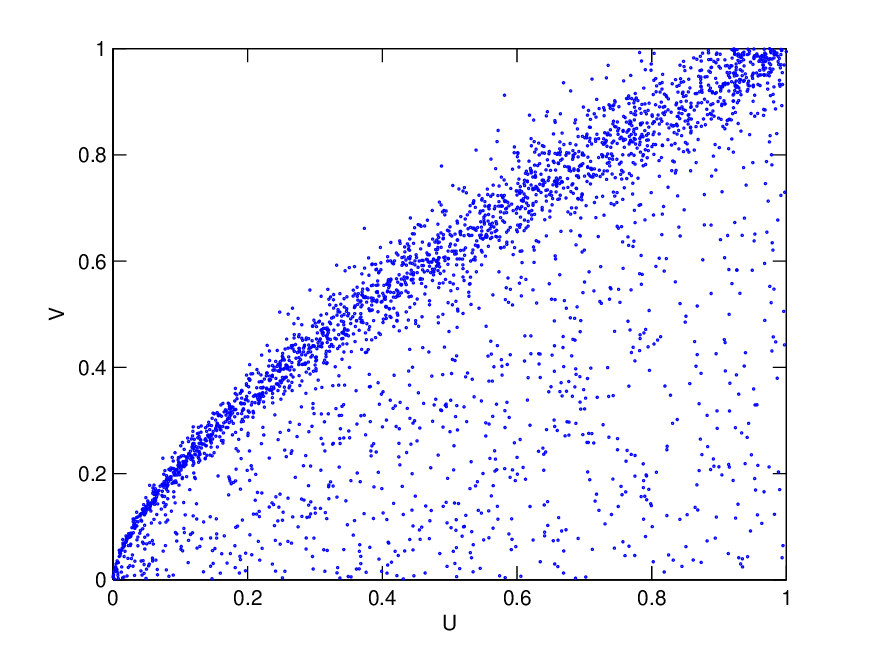}
\hfill
\includegraphics[width=7.5cm,height=7.5cm]{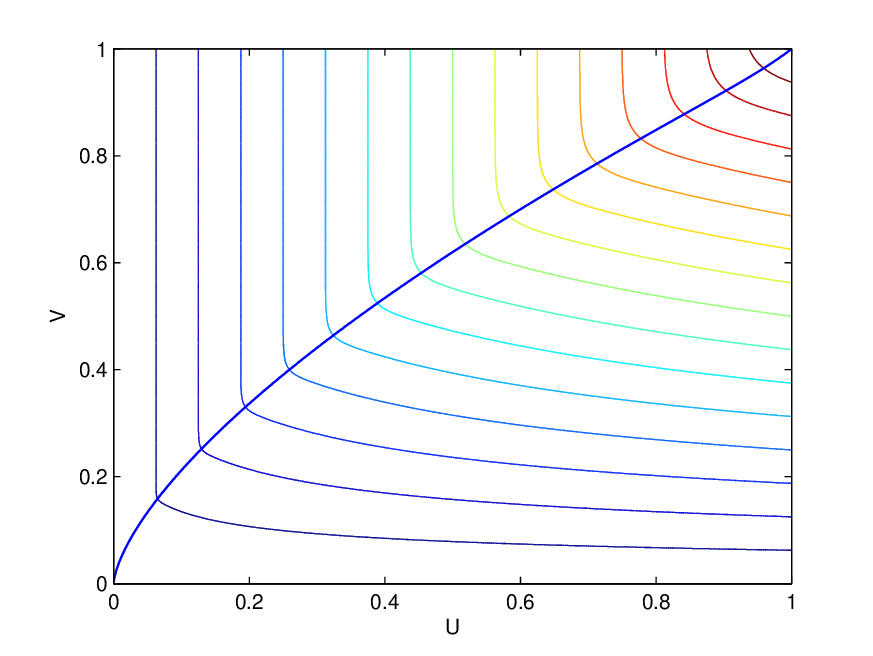}
\caption{Simulated generalized Clayton copula (left) and its contours (right) when $\gamma_0=0.04$ and $\gamma_1=0.02$, with the path of maximal dependence $(\varphi^*(u),u^2/\varphi^*(u))_{0\le u \le 1}$ superimposed on the right-hand panel.}
\label{fig-genclayton}
\end{figure}

The generalized Clayton copula is given by the formula
\[
C(u,v)=u^{\gamma_1/\tilde{\gamma}_1}(u^{-1/\tilde{\gamma}_1}
+v^{-1/{\gamma}_0}-1)^{-\gamma_0}
\quad \textrm{for} \quad 0\le u,v\le 1,
\]
where $\gamma_0>0$ and $\gamma_1\ge 0$ are parameters, and $\tilde{\gamma}_1=\gamma_0+\gamma_1$. Tedious computations, which we have relegated to Appendix \ref{appendix}, show that the function of maximal dependence $\varphi^*$ (Figure \ref{fig-genclayton})
is unique and satisfies the equation
\begin{equation}
\zeta(\varphi^*(u))=0
 \label{qq-3ef}
\end{equation}
for every $0\le u \le 1$, where
\[
\zeta(x)=x^{-1/{\gamma}_0}\left(x^{-1/\tilde{\gamma}_1}-(\gamma_1 / \tilde{\gamma}_1) \right)-\left(1-(\gamma_1 / \tilde{\gamma}_1) \right)u^{-2/{\gamma}_0}
\]
for all $x\in [u^2,1]$. Further tedious calculations (Appendix \ref{appendix}) show that the lower-tail index of maximal dependence is
\begin{equation}
\kappa^*_L =1+\frac{\gamma_1}{\gamma_1+2\gamma_0}.
\label{qq-3e}
\end{equation}

\subsection{Archimedean copula}

Another example illustrating the difficulty of deriving closed-form expressions for paths of maximal dependence is provided by the Archimedean copula
\begin{equation}
C(u,v)=\psi^{-1}(\psi(u)+\psi(v))
\quad \textrm{for} \quad 0\le u,v\le 1,
\label{arch-copula}
\end{equation}
where $\psi:[0,1] \rightarrow [0,\infty]$ is such that $\psi(1)=0$, $\psi^{'}(u)<0$, and $\psi^{''}(u)>0$ for all $0<u<1$. Assume that the Archimedean generator $\psi $ is strict, that is, $\psi(u)\to \infty$ when $u\downarrow 0$, as only in this case can we have a path of maximal dependence. For every $u\in [0,1]$, to determine those $x\in [u^2,1]$ for which the equation $(\partial /\partial x )C(x, u^2/x) = 0$ holds, we notice that when the function $x\psi'(x)$ is increasing on $[u^2,1]$, which we assume from now on, then
\[
x\psi'(x)-(u^2/x)\psi'(u^2/x)= ~\begin{cases}
<0 & \textrm{ when}\quad x< u^2/x, \\
>0 & \textrm{ when}\quad x>u^2/x. \\
\end{cases}
\]
Thus, $\psi(x)+\psi(u^2/x) \geq 2\psi(u)$ for all $x\in[u^2,1]$ or, equivalently,
\[
\psi^{-1}(\psi(x)+\psi(u^2/x)) \leq \psi^{-1}(2\psi(u)).
\]
Therefore, when the function $x\psi'(x)$ is  increasing on $[u^2,1]$, then the diagonal is a path of maximal dependence. In this case, the Archimedean copula is PQD. We observe that all PQD Archimedean copulas documented by Nelsen (2006) have increasing functions $x\psi'(x)$.

\section{Concluding notes}
\label{section-7}

We have demonstrated that the classical indices of tail dependence, which are based on the behaviour of copulas along the diagonal path, do not generally capture the maximal degree of tail dependence. For this reason, we have proposed conservative indices of tail dependence that hinge on the new notion of \textit{paths of maximal dependence} that we have introduced herein. We have used a number of specific copulas, as well as a numerical example, to elucidate relevant main ideas both analytically and numerically. Our approach to assessing the tail dependence in copulas conforms with the new paradigm of prudence in modern quantitative risk management (cf.\, ORSA, 2014).

A number of challenging problems remain in this area. For example, statistical inference for the indices would be of interest, and we believe that tools developed in the area of M-estimators and related weak convergence theorems would be helpful (cf., e.g., van der Vaart and Wellner, 1996; van de Geer, 2009; and references therein). Another interesting problem would be to explore uniqueness and other properties of the path(s) of maximal dependence, and we think that results discussed by Dharmadhikari and Joag-dev (1987) would be particularly helpful.

\section*{Acknowledgments}

We are grateful to two anonymous referees for constructive criticism, suggestions, and thoughts provoking questions, as well as to Prof.~Dr.~Paul Embrechts and all participants of the ETH's Series of Talks in Financial and Insurance Mathematics for feedback and insights.

Our research has been supported by the Natural Sciences and Engineering Research Council (NSERC) of Canada. Jianxi Su also acknowledges the financial support of the Government of Ontario via its Ontario Graduate Scholarship program.

We enjoyed our respective working visits at York University and the University of Western Ontario, and we thank both institutions for much appreciated research environments and hospitality.

\section*{References}
\def\hang{\hangindent=\parindent\noindent}

\hang
\textsc{Asimit, A.V. and Badescu, A.L.} (2010) Extremes on the discounted aggregate claims in a time dependent risk model. \textit{Scandinavian Actuarial Journal} 2010, 93--104.

\hang
\textsc{Asimit, A.V., Furman, E. and Vernic, R.} (2010)
On a multivariate Pareto distribution.
\textit{Insurance: Mathematics and Economics}
46(2), 308--316.

\hang
\textsc{Asimit, A.V., Furman, E., Tang, Q. and Vernic, R.} (2011) Asymptotics for risk capital allocations based on conditional tail expectation. \textit{Insurance: Mathematics and Economics} 49, 310--324.

\hang
\textsc{Asimit, A.V. and Jones, B.L.} (2008) Dependence and the asymptotic behavior of large claims reinsurance. \textit{Insurance: Mathematics and Economics} 43, 407--411.

\hang
\textsc{Balakrishnan, N. and Lai, C.D.} (2009) \textit{Continuous Bivariate Distributions.} (Second Edition.) Springer. New York.

\hang
\textsc{Biener, C. and Eling, M.} (2012) Insurability in microinsurance markets: an analysis of problems and potential solutions. \textit{The Geneva Papers} 37, 77--107.

\hang
\textsc{Churchill, C. and Matul, M.} (2012) \textit{Protecting the Poor: A Microinsurance Compendium.} (Volume II.) International Labour Organization. Geneva, Switzerland.

\hang
\textsc{Coles, S., Heffernan, J. and Tawn, J.} (1999) Dependence measures for extreme value analyses. \textit{Extremes} 2, 339--365.

\hang
\textsc{Dharmadhikari, S. and Joag-dev, K.} (1987) \textit{Unimodality, Convexity and Applications.} Academic press, New York.

\hang
\textsc{Durante, F., Fern\'{a}ndez-S\'{a}nchez, J. and Pappad\'{a}, R.} (2014) Copulas, diagonals, and tail dependence. \textit{Fuzzy Sets and Systems} (in print).

\hang
\textsc{Durante, F., Pappad\`{a}, R. and Torelli, N.} (2014) Clustering of financial time series in risky scenarios.  \textit{Advances in Data Analysis and Classification} 8, 359--376.

\hang
\textsc{Embrechts, P.} (2009) Copulas: a personal view. \textit{Journal of Risk and Insurance} 76, 639--650.

\hang
\textsc{Embrechts, P., Lindskog, F. and McNeil, A.} (2003) Modeling dependence with copulas and applications to risk management. In: \textit{Handbook of Heavy Tailed Distributions in Finance} (ed. S. Rachev) Chapter 8, pp. 329--384. Elsevier.

\hang
\textsc{Embrechts, P., Puccetti, G. and R\"{u}schendorf, L.} (2013) Model uncertainty and VaR aggregation. \textit{Journal of Banking and Finance} 37, 2750--2764.

\hang
\textsc{Embrechts, P., Puccetti, G., R\"{u}schendorf, L., Wang, R. and  Beleraj, A.} (2014) An academic response to Basel 3.5. \textit{Risks} 2, 25--48.

\hang
\textsc{Fischer, M.J. and Klein, I.} (2007) Some results on weak and strong tail dependence coefficients for means of copulas.
Diskussionspapiere // Friedrich-Alexander-Universit\"{a}t Erlangen-N\"{u}rnberg, Lehrstuhl f\"{u}r
Statistik und \"{O}konometrie, No. 78/2007.

\hang
\textsc{Frees, E.W.} (2010) \textit{Regression Modeling with Actuarial and Financial Applications.} Cambridge University Press, New York.

\hang
\textsc{Frees, E.W. and Wang, P.} (2005) Credibility using copulas. \textit{North American Actuarial Journal} 9(2), 31--48.

\hang
\textsc{Frees, E.W., Carriere, J. and Valdez, E.} (1996) Annuity valuation with dependent mortality. \textit{Journal of Risk and Insurance} 63, 229--261.

\hang
\textsc{Furman, E., Kuznetsov, A., Su, J. and Zitikis, R.} (2015) \textit{The Maximal Tail Dependence Path of the Gaussian Copula is Diagonal.} Technical report available at SSRN: http://ssrn.com/abstract=2558675

\hang
\textsc{Furman, E. and Zitikis, R.} (2010)
General Stein-type covariance decompositions with applications to insurance and finance.
\textit{ASTIN Bulletin} 40, 369--375.

\hang
\textsc{Heffernan, J.E.} (2000) A directory of coefficients of tail dependence. \textit{Extremes} 3, 279--290.

\hang
\textsc{Hua, L. and Joe, H.} (2014) Strength of tail dependence based on conditional tail expectation. \textit{Journal of Multivariate Analysis} 123, 143--159.

\hang
\textsc{Jaworski, P., Durante, F. and H\"{a}rdle, W.} (2013) \textit{Copulae in Mathematical and Quantitative Finance.} Springer. Berlin.

\hang
\textsc{Jaworski, P., Durante, F., H\"{a}rdle, W. and Rychlik, T.} (2010) \textit{Copulas Theory and its Applications.} Springer. Berlin.

\hang
\textsc{Joe, H., Li, H. and Nikoloulopoulos A.K.} (2010) Tail dependence functions and vine copulas. \textit{Journal of Multivariate Analysis} 101, 252--270.

\hang
\textsc{Ledford, A.W. and Tawn, J.A.} (1996) Statistics for near independence in multivariate extreme values. \textit{Biometrika} 83, 169--187.

\hang
\textsc{Li, H. and Wu, P.} (2013) Extremal dependence of copulas: a tail density approach. \textit{Journal of Multivariate Analysis} 144, 99--111.

\hang
\textsc{McNeil, A. J., Frey, R. and Embrechts, P.} (2005) \textit{Quantitative Risk Management.} Princeton University Press. Princeton.

\hang
\textsc{Nelsen, R.B.} (2006)
\textit{An Introduction to Copulas.} (Second Edition.) New York: Springer.

\hang
\textsc{Nikoloulopoulos, A.K. and Karlis, D.} (2008) Multivariate logit copula model with an application to dental data. \textit{Statistics in Medicine} 27, 6393--6406.

\hang
ORSA (2014). \textit{Own Risk and Solvency Assessment.} Office of the Superintendent of Financial Institutions, Government of Canada, Ottawa. http://www.osfi-bsif.gc.ca/eng/fi-if/rg-ro/gdn-ort/gl-ld/Pages/e19.aspx.  Accessed on March 17, 2015.

\hang
\textsc{Patton, A.J.} (2012) A review of copula models for economic time series.  \textit{Journal of Multivariate Analysis} 110, 4--18.

\hang
\textsc{Puccetti, G. and R\"{u}schendorf, L.} (2014) Asymptotic equivalence of conservative value-at-risk- and expected shortfall-based capital charges. \textit{Journal of Risk} 16, 3--22.

\hang
\textsc{van der Vaart, A.W. and Wellner, J.A.} (1996) \textit{Weak Convergence and Empirical Processes. With Applications to Statistics.} Springer, New York.

\hang
\textsc{van de Geer, S.A.} (2009) \textit{Empirical Processes in M-Estimation.} Cambridge University Press, Cambridge.

\hang
\textsc{Weng, C. and Zhang, Y.} (2012) Characterization of multivariate heavy-tailed distribution families via copula. \textit{Journal of Multivariate Analysis} 106, 178--186.

\appendix
\section{Appendix: proofs}
\label{appendix}

\begin{proof}[Proof of Theorem \ref{cont}]
We prove by contradiction. Hence, assume that there is a point $u_0$ of  discontinuity of the function $\varphi^*$ of maximal dependence. This implies that there is a sequence $u_n\to u_0$ such that $\varphi^*(u_n) \to x_0 \not= \varphi^*(u_0)$. By the continuity of the function $(v,w) \mapsto C(v,w)$ (cf., e.g., Embrechts et al, 2003; Nelsen, 2006), we have
\begin{equation}\label{cc-1}
C(\varphi^*(u_n),u_n^2/\varphi^*(u_n))\to C(x_0,u_0^2/x_0).
\end{equation}
By the uniqueness of the function of maximal dependence, $x_0$ cannot maximize the function $x\mapsto C(x,u_0^2/x)$ because otherwise $x_0$ would be equal to $\varphi^*(u_0)$. On the other hand, we have
\begin{equation}\label{cc-2}
C(\varphi^*(u_n),u_n^2/\varphi^*(u_n))
=\max_{x\in [u_n^2,1]}C(x,u_n^2/x)
\to \max_{x\in [u_0^2,1]}C(x,u_0^2/x),
\end{equation}
where the limit holds due to the following reasons. When $u_n<u_0$, then
\begin{align*}
\max_{x\in [u_n^2,1]}C(x,u_n^2/x)
&=\max \Big \{ \max_{x\in [u_n^2,u_0^2)}C(x,u_n^2/x),
\max_{x\in[u_0^2,1]}C(x,u_n^2/x) \Big \}
\\
&\to\max \Big \{ u_0^2,
\max_{x\in[u_0^2,1]}C(x,u_0^2/x) \Big \}
\\
&= \max_{x\in [u_0^2,1]}C(x,u_0^2/x),
\end{align*}
{\color{blue}where the right-most equation holds due to $u_0^2 = C(x,u_0^2/x)$ for $x=1$.} When $u_n>u_0$, then
\begin{align*}
&\lim_{n\to \infty }\Big | \max_{x\in [u_n^2,1]}C(x,u_n^2/x)
-\max \Big \{ \max_{x\in [u_0^2,u_n^2)}C(x,u_0^2/x),
\max_{x\in[u_n^2,1]}C(x,u_0^2/x) \Big \} \Big |
\\
&=\lim_{n\to \infty }\Big | \max_{x\in [u_n^2,1]}C(x,u_n^2/x)
-\max \Big \{ u_0^2,
\max_{x\in[u_n^2,1]}C(x,u_0^2/x) \Big \} \Big |
\\
&=\lim_{n\to \infty }\Big | \max_{x\in [u_n^2,1]}C(x,u_n^2/x)
-\max_{x\in[u_n^2,1]}C(x,u_0^2/x)  \Big |
\\
&\le \lim_{n\to \infty }\max_{x\in [u_n^2,1]}\big | C(x,u_n^2/x)
-C(x,u_0^2/x)  \big |
=0
\end{align*}
because of the uniform continuity of the function $(v,w) \mapsto C(v,w)$ (cf., e.g., Embrechts et al, 2003; Nelsen, 2006). This completes the proof of limit (\ref{cc-2}). Combining statements (\ref{cc-1}) and (\ref{cc-2}), we get the equation
\begin{equation}\label{cc-3}
\max_{x\in [u_0^2,1]}C(x,u_0^2/x)=C(x_0,u_0^2/x_0),
\end{equation}
which implies that $x_0$ must be equal to $\varphi^*(u_0)$ because of the uniqueness of the maximum. We have arrived at a contradiction, which establishes the continuity of $\varphi^*$ and completes the proof of Theorem \ref{cont}.
\end{proof}

Throughout the remaining proofs, though perhaps not mentioned explicitly, $u$ and $v$ are always in the interval $[0,1]$, and $x$ is always in the interval $[u^2,1]$.

\begin{proof}[Proof of equations (\ref{mo-1a0}) and (\ref{qq-3eee})]
Without loss of generality, let $a<b$. Then $x_{a,b}=u^{2b/(a+b)}$ is smaller than $x_{b,a}=u^{2a/(a+b)}$, and so the two functions $x\mapsto C_{a,b}(x,u^2/x)$ and $x\mapsto C_{b,a}(x,u^2/x)$ are increasing on the interval $(u^2,x_{a,b})$ and decreasing on $(x_{b,a},1)$. Hence, the maximum of the function $x\mapsto C(x,u^2/x)$ can only be achieved on the interval $[x_{a,b},x_{b,a}]$, where we have the formula
\[
C(x,u^2/x)={1\over 2} \Big (  u^{2(1-a)}x^a+u^2 / x^a \Big ).
\]
We split the interval $[x_{a,b},x_{b,a}]$ into two subintervals:  $[x_{a,b},u]$ where the function $x\mapsto C(x,u^2/x)$ is decreasing, and  $[u,x_{b,a}]$ where the function is increasing. From this we conclude that the function $x\mapsto C(x,u^2/x)$ can achieve its maximum only at  $x_{a,b}$ and/or $x_{b,a}$, and it attains the same value at both end-points. Consequently, the two end-points are maxima and thus define two functions of maximal dependence, which are $\varphi_1^*$ and $\varphi_2^*$ given by formulas (\ref{mo-1a0}). For any of these two functions, $\varphi_k^*$, we have
\[
C(\varphi_k^*(u),u^2/\varphi_k^*(u))={1\over 2} \Big ( u^{2-2(a\wedge b)^2/(a+b)}+u^{2-2ab/(a+b)} \Big ).
\]
With $\kappa =2- 2ab/(a+b)$ we check that
\[
{\Pi^*(u) \over u^{\kappa } }=
{C(\varphi_k^*(u),u^2/\varphi_k^*(u)) \over u^{\kappa } }\to {1\over 2}
\]
when $u \downarrow 0 $. Consequently, the lower-tail index of maximal dependence $\kappa^*_L$ is equal to $\kappa $. This concludes the proof of equations (\ref{mo-1a0}) and (\ref{qq-3eee}).
\end{proof}

\begin{proof}[Proof of equations (\ref{qq-3ef}) and (\ref{qq-3e})]
We start with the equation
\begin{equation}
C(x,u^2/x)=x^{\gamma_1/\tilde{\gamma}_1}
(x^{-1/\tilde{\gamma}_1}+u^{-2/\gamma_0}x^{1/\gamma_0}-1)^{-\gamma_0}.
\label{qq-1}
\end{equation}
Finding the maximum of $C(x,u^2/x)$ with respect to $x\in [u^2,1]$ is the same as finding the maximum of its logarithm $\log C(x,u^2/x)$, which is easier in this particular case. We have
\[
\frac{\partial}{\partial x}\log(C(x,u^2/x)) =\frac{x^{-1/\tilde{\gamma}_1}
+((\gamma_1 / \tilde{\gamma}_1)-1)(u^2/x)^{-1/{\gamma}_0}-(\gamma_1 / \tilde{\gamma}_1)}{x(x^{-1/\tilde{\gamma}_1}+(u^2/x)^{-1/{\gamma}_0}-1)}.
\]
The denominator $x(x^{-1/\tilde{\gamma}_1}+(u^2/x)^{-1/{\gamma}_0}-1)$ is positive for all $x\in [u^2, 1]$, and thus we need to find those $x\in [u^2, 1]$ that make the numerator equal to $0$. This is equivalent to solving the equation
$\zeta(x)=0$. The solution to this equation is unique because
\[
\zeta(u^2)=u^{-2/{\gamma}_0}(u^{-2/\tilde{\gamma}_1}-1)>0,
\]
\[
\zeta(1)=\left(1-(\gamma_1 / \tilde{\gamma}_1)\right)-\left(1-(\gamma_1 / \tilde{\gamma}_1)\right)u^{-2/\tilde{\gamma}_1}<0,
\]
and, for all $x\in [u^2, 1]$,
\[
\zeta'(x)=-\frac{2}{\tilde{\gamma}_1}x^{-1/\tilde{\gamma}_1-1/{\gamma}_0-1} - \frac{\gamma_1}{\tilde{\gamma}_1 \gamma_0} \left(x^{-1/\tilde{\gamma}_1-1/{\gamma}_0-1} -x^{-1/\gamma_0-1}\right)<0.
\]

We cannot derive a closed form solution to the equation $\zeta(x)=0$, but we already know that the solution $x=\varphi^*(u)\in [u^2,1]$ exists and is unique. Furthermore, the solution satisfies the equation
\begin{equation}
\label{easy-form}
x^{-1/\tilde{\gamma}_1-1/\gamma_0}\left(1-(\gamma_1 / \tilde{\gamma}_1)x^{1/\tilde{\gamma}_1} \right) =\left(1- (\gamma_1 / \tilde{\gamma}_1)\right)u^{-2/\gamma_0}.
\end{equation}
We also have that $1-(\gamma_1 / \tilde{\gamma}_1)x^{1/\tilde{\gamma}_1} \in [1-(\gamma_1 / \tilde{\gamma}_1),1]$. From these facts we conclude that $x=\varphi^*(u)\downarrow 0$ when $u \downarrow 0$.

Due to the lack of closed-form expression for the function of maximal dependence, we cannot obtain a closed-form expression for the maximal  probability $\Pi^*(u)$ either. Nevertheless, we can obtain a closed-form expression for $\kappa^*_L$. Starting with equation (\ref{easy-form}), we arrive at the following one:
\begin{equation}
\label{easy-form-1}
x= \bigg ( u^{2/\gamma_0}\left(1-(\gamma_1 / \tilde{\gamma}_1)x^{1/\tilde{\gamma}_1} \right)\Big / \left(1- (\gamma_1 / \tilde{\gamma}_1)\right)
\bigg )^{\tilde{\gamma}_1\gamma_0/(\tilde{\gamma}_1+\gamma_0)}.
\end{equation}
Denote $r(x)=\left(1-(\gamma_1 / \tilde{\gamma}_1)x^{1/\tilde{\gamma}_1}\right) /\left(1- (\gamma_1 / \tilde{\gamma}_1)\right)$. Replacing all the $x$'s on the right-hand side of equation (\ref{qq-1}) by the right-hand side of equation (\ref{easy-form-1}), we obtain
\begin{align*}
C(x,u^2/x)&=\big \{u^{2/\gamma_0}r(x)
\big \}^{\tilde{\gamma}_1\gamma_0/(\tilde{\gamma}_1+\gamma_0)}
\\
 &=u^{\frac{2\gamma_1}{\tilde{\gamma}_1+\gamma_0}}
 r(x)^{\frac{\gamma_0\gamma_1}{\tilde{\gamma}_1+\gamma_0}}
 \bigg(r(x)^{-\frac{\gamma_0}{\tilde{\gamma}_1+\gamma_0}}u^{-\frac{2}{\tilde{\gamma}_1+\gamma_0}}+r(x)^{\frac{\tilde{\gamma}_1/\gamma_0}{\tilde{\gamma}_1+\gamma_0}} u^{-\frac{2}{\tilde{\gamma}_1+\gamma_0}}-1 \bigg)^{-\gamma_0}.
\end{align*}
Consequently, with
$\kappa =1+\gamma_1/(\tilde{\gamma}_1+\gamma_0)$
we have
\[
{C(x,u^2/x) \over u^{\kappa } }
=r(x)^{\frac{\gamma_0\gamma_1}{\tilde{\gamma}_1+\gamma_0}}\bigg ( r(x)^{-\frac{\gamma_0}{\tilde{\gamma}_1+\gamma_0}}+r(x)^{\frac{\tilde{\gamma}_1/\gamma_0}{\tilde{\gamma}_1+\gamma_0}} -u^{\frac{2}{\tilde{\gamma}_1+\gamma_0}} \bigg)^{-\gamma_0} \to c \in (0,\infty)
\]
when $u \downarrow 0 $. This proves that the lower-tail index of maximal dependence $\kappa^*_L$ is equal to $\kappa $. The proof of equations (\ref{qq-3ef}) and (\ref{qq-3e}) is finished.
\end{proof}

\end{document}